\documentclass[a4paper,12pt]{article}
\usepackage{amssymb}
\usepackage{amsthm}
\newtheorem{Lemma}{Lemma}
\newtheorem{Proposition}{Proposition}
\newtheorem{Theorem}{Theorem}

\newtheorem{Corollary}{Corollary}

\newcommand{\bP}{\mathbb{P}}
\newcommand{\Q}{\mathbb{Q}}

\begin{document}
\title{On Waring's problem \\ for partially symmetric tensors \\
\emph{(variations on a theme of Mella)}}
\author{Claudio~Fontanari}
\date{}
\maketitle 
\begin{small}
\begin{center}
\textbf{Abstract}
\end{center}
Here we address a Waring type problem for partially symmetric tensors, 
extending previous work by Massimiliano Mella in the totally symmetric 
case of forms. In particular, we provide an explicit answer in lower 
dimensional cases. 
\vspace{0.5cm}

\noindent
\textsc{A.M.S. Math. Subject Classification (2000)}: 14N05.

\end{small}
\vspace{0.5cm}

\section{Introduction}
Waring's problem for forms, concerning the decomposition of a generic 
homogeneous polynomial as a sum of powers of linear forms, has now a 
complete and rigorous solution thanks to the spectacular results
obtained by by Alexander and Hirschowitz in \cite{AleHir:95}. 
However, several related problems are still widely open: 
for instance, one would like to determine all cases in which the 
above decomposition is unique. As far as we know, the best achievement in 
this direction is the following Theorem, recently proved by Massimiliano 
Mella (see \cite{Mella:04}, Theorem~1):

\begin{Theorem}\label{mella} \emph{\textbf{(Mella)}} 
Fix integers $d > r > 1$ and $k \ge 1$ such that 
$(k+1)(r+1) = {{r+d} \choose r}$. 
Then the generic homogeneous polynomial of degree $d$ in $r+1$ variables 
can be expressed as a sum of $k+1$ $d$-th powers of linear forms in a 
unique way if and only if $d=5$ and $r=2$.   
\end{Theorem}

It seems rather interesting also to consider variants of the original 
Waring's problem, the best known one being the so-called Extended Waring's 
Problem (see \cite{Ciliberto:01}, Problem~7.6, and \cite{Fontanari:02}), 
which asks for a simultaneous decomposition of several forms. Here 
instead we propose another generalization, which regards forms as 
symmetric tensors and is concerned with the decomposition of a 
wider class of tensors. More precisely, for fixed integers $n \ge 1$, 
$r \ge 1$, and $\underline{d} = (d_1, \ldots, d_n)$, let $V_{r, d_i}$ 
be the Veronese embedding of $\bP^r$ of degree $d_i$ and let $\Sigma_{r, 
\underline{d}}$ denote the Segre embedding of $V_{r, d_1} \times 
\ldots \times V_{r, d_n}$ into $\bP^N$, $N = \prod_{i=1}^n 
{{r+d_i} \choose r} - 1$. According to \cite{CGG:03}, \S~4, we can 
view $\bP^N$ inside $\bP^M$, $M = (r+1)^{d_1+ \ldots + d_n}-1$, as 
the space of tensors which are invariant with respect to the natural 
actions of the symmetric groups $S_{d_1}, \ldots, S_{d_n}$ on the 
coordinates of $\bP^M$. Therefore we identify the points of $\bP^N$ 
with the $(d_1, \ldots, d_n)$ partially symmetric tensors in $\bP^M$; 
in particular, the points on the Segre-Veronese variety $\Sigma_{r, 
\underline{d}}$ correspond to the partially symmetric tensors in 
$\bP^M$ which are decomposable. Now we can state a general result
about decomposition of partially symmetric tensors 
(we refer to \cite{ChiCil:02} for the definitions of defectivity and 
weak defectivity):  

\begin{Theorem}\label{main}  
Fix integers $n \ge 1$ and $r \ge 1$ with $n r \ge 2$, 
$d_n \ge d_{n-1}\ge \ldots \ge d_1 \ge r+1$ 
and $k \ge 1$ such that 
$$
\prod_{i=1}^n {{r+d_i} \choose r} = (n r + 1)(k+1).
$$
If $\Sigma_{r, \underline{d}}$ is neither $k$-defective nor
$(k-1)$-weakly defective, then a general $(d_1, \ldots, d_n)$ 
partially symmetric tensor in $\bP^{(r+1)^{d_1+ \ldots +d_n}-1}$
can be expressed in $\nu \ge 2$ ways as a sum of $k+1$ decomposable 
tensors.
\end{Theorem}

Notice that the special case $n=1$ of Theorem~\ref{main} is the main 
ingredient of Mella's Theorem~\ref{mella} (some additional work is 
needed in order to check the assumption about weak defectivity, see 
\cite{Mella:04}, \S~3 and \S~4). 
Next we point out a couple of new results, obtained from 
Theorem~\ref{main} together with previous contributions 
(see \cite{ChiCil:02} and \cite{CGG:03}). 
The easiest one is the following: 

\begin{Corollary}\label{two}
Fix integers $d_2 \ge d_1 \ge 4$ and $k \ge 1$ such that 
$(d_1+1)(d_2+1)=3(k+1)$. Then a general $(d_1,d_2)$ partially 
symmetric tensor in $\bP^{2^{d_1+d_2}-1}$ can be expressed in 
$\nu \ge 2$ ways as a sum of $k+1$ decomposable tensors. 
\end{Corollary}

A careful application of the so-called Horace method (see 
Proposition~\ref{horace} and Proposition~\ref{weakly}) yields 
also the following:

\begin{Corollary}\label{three}
Fix integers $d_3 \ge d_2 \ge d_1 \ge 3$ and $k \ge 1$ such that 
$(d_1+1)(d_2+1)=4(k+1)$. Assume that 
\begin{equation}\label{assumption} 
k+1 \le (d_3-2) \left[ \frac{(d_1+1)(d_2+1)}{3} \right]
\end{equation}
where $\left[ . \right]$ denotes the integral part.
Then a general $(d_1,d_2,d_3)$ partially 
symmetric tensor in $\bP^{2^{d_1+d_2+d_3}-1}$ can be expressed in 
$\nu \ge 2$ ways as a sum of $k+1$ decomposable tensors. 
\end{Corollary}

We work over an algebraically closed field of characteristic zero.

This research is part of the T.A.S.C.A. project of I.N.d.A.M., 
supported by P.A.T. (Trento) and M.I.U.R. (Italy).

\section{The proofs}
The proof of Theorem~\ref{main} involves two main ingredients. 

The first result is well-known to the specialists (indeed,  
Bronowski was aware of it already in 1932: see \cite{Bronowski:32} 
for the original statement and \cite{CMR:04}, Corollary~4.2,
\cite{CR:04}, Theorem~2.7, \cite{Mella:04}, Theorem~2.1, 
for rigorous modern proofs):

\begin{Lemma}\label{bronowski}
Let $X \subset \bP^{r(k+1)+k}$ be a smooth and irreducible projective 
variety of dimension $r$ such that through the general point of 
$\bP^{r(k+1)+k}$ there is exactly one $\bP^k$ which is $(k+1)$-secant 
to $X$. Then the projection of $X$ to $\bP^r$ from a general tangent
space to the $(k-1)$-secant variety $S^{k-1}(X)$ is birational. 
\end{Lemma} 

The second result is a version of the so-called Noether-Fano inequality 
for Mori fiber spaces (see \cite{Corti:95}, Definition~(3.1) and 
Theorem~(4.2)):

\begin{Lemma}\label{noether-fano}
Let $\pi: X \to S$ and $\rho: Y \to T$ be two Mori fiber spaces and 
let $\varphi: X \dashrightarrow Y$ be a birational not biregular map. 
Choose a very ample linear system $\mathcal{H}_Y$ on $Y$, 
define $\mathcal{H}_X = \varphi^*(\mathcal{H}_Y)$ and let 
$\mu \in \Q$ such that $\mathcal{H}_X \equiv - \mu K_X + \varphi^*(A)$ 
for some divisor $A$ on $S$. Then either $(X, 1/\mu \mathcal{H}_X)$ 
has not canonical singularities or $K_X + 1/\mu \mathcal{H}_X$ is 
not nef.
\end{Lemma}

\emph{Proof of Theorem~\ref{main}.} Since $\Sigma_{r, \underline{d}}$ 
lies in a projective space $\bP^N$ of dimension equal to the expected 
dimension of the $k$-secant variety $S^k(\Sigma_{r, \underline{d}})$, 
the non $k$-defectivity of $\Sigma_{r, \underline{d}}$ implies that 
through the general point of $\bP^N$ there is at least one $\bP^k$ 
which is $(k+1)$-secant to $\Sigma_{r, \underline{d}}$. Hence the 
existence of at least one expression as in the statement follows from 
the translation described in the Introduction and we have only to check 
that such an expression is never unique. Arguing by contradiction, 
assume that through the general point of $\bP^N$ there is exactly one 
$\bP^k$ which is $(k+1)$-secant to $\Sigma_{r, \underline{d}}$. 
From Lemma~\ref{bronowski} we obtain a birational map
$$
\bP^r \times \ldots \times \bP^r \dashrightarrow \bP^{n r}  
$$
induced by the linear system $\mathcal{H} = \mathcal{O}_{\bP^r}(d_1) 
\otimes \ldots \otimes \mathcal{O}_{\bP^r}(d_n)(-2 p_1 \ldots -2 p_k)$, 
where the $p_i$'s are general points in $\bP^r \times \ldots \times \bP^r$. 
Since $\Sigma_{r, \underline{d}}$ is not $(k-1)$-weakly defective, 
by \cite{ChiCil:02}, Theorem~1.4, a general divisor $H \in \mathcal{H}$ 
has only ordinary double points at $p_1, \ldots, p_k$ and is elsewhere 
smooth. In particular, $H$ is irreducible and after the blow-up of 
$p_1, \ldots, p_k$ it becomes smooth, hence we see that $H$ has only 
canonical singularities. Since both $\bP^r \times \ldots 
\times \bP^r$ and $\bP^{n r}$ regarded as projective bundles are Mori 
fiber spaces, we can apply Lemma~\ref{noether-fano} with $\mu := 
\frac{d_1}{r+1}$. By our numerical assumptions, we have 
$-(r+1)+\frac{1}{\mu} d_i \ge 0$ for every $i$, hence $K + 1/\mu 
\mathcal{H}$  is nef and this contradiction ends the proof. 

\qed 

\emph{Proof of Corollary~\ref{two}.} According to Theorem~\ref{main}, 
we have only to show that the Segre embedding of $V_{1, d_1} \times 
V_{1, d_2}$ is neither $k$-defective nor $(k-1)$-defective. 
These two properties can be easily checked by looking at the 
Classification Theorem~1.3 of \cite{ChiCil:02}, so the proof 
is over.

\qed

\begin{Proposition}\label{horace}
Fix integers $n \ge 1$, $r_1, \ldots, r_n \ge 1$, $d_1, \ldots, d_n \ge 1$,
$d_{n+1} \ge 2$, 
$0 \le h \le \left[ \frac{ \prod_{i=1}^n {{r_i + d_i} \choose r_i}}
{(r_1 + \ldots + r_n+1)} \right]$, 
$0 \le l \le \left[ \frac{ \prod_{i=1}^{n} {{r_i + d_i} \choose r_i}
(d_{n+1}+1)}{(r_1 + \ldots + r_n + 2)} \right]
- \left[ \frac{ \prod_{i=1}^n {{r_i + d_i} \choose r_i}}
{(r_1 + \ldots + r_n+1)} \right]$. 
Choose $l$ general points $p_1, \ldots, p_l$ in 
$\bP^{r_1} \times \ldots \times \bP^{r_n} \times \bP^1$,
let $D$ be a divisor of type $(0, \ldots, 0, 1)$ on 
$\bP^{r_1} \times \ldots \times \bP^{r_n} \times \bP^1$, and pick $h$ 
general points $q_1, \ldots, q_h$ on $D$. Assume that  
\begin{eqnarray*}
\dim \vert \mathcal{O}_{\bP^{r_1} \times \ldots \times \bP^{r_n} \times \bP^1}
(d_1, \ldots, d_n, d_{n+1}-1) \left(-\sum_{i=1}^l 2 p_i\right) \vert = \\ 
= \prod_{i=1}^n {{r_i + d_i} \choose r_i}(d_{n+1})
- (r_1 + \ldots + r_n + 2)l
\end{eqnarray*}
and 
\begin{eqnarray*}
\dim \vert \mathcal{O}_{\bP^{r_1} \times \ldots \times \bP^{r_n}}
(d_1, \ldots, d_n) \left(-\sum_{j=1}^h 2 q_j\right) \vert = \\
= \prod_{i=1}^n {{r_i + d_i} \choose r_i} - (r_1 + \ldots + r_n + 1)h,
\end{eqnarray*}
as expected, and that 
\begin{eqnarray*}
\dim \vert \mathcal{O}_{\bP^{r_1} \times \ldots \times \bP^{r_n} \times \bP^1}
(d_1, \ldots, d_n, d_{n+1}-2)\left(-\sum_{i=1}^l 2 p_i\right) \vert  
\le \prod_{i=1}^n {{r_i + d_i} \choose r_i}(d_{n+1}) \\
- (r_1 + \ldots + r_n + 2)l - h.
\end{eqnarray*}
Then we have
\begin{eqnarray*}
\dim \vert \mathcal{O}_{\bP^{r_1} \times \ldots \times 
\bP^{r_n} \times \bP^1}(d_1, \ldots, d_n, d_{n+1}) \left(-\sum_{i=1}^l 2 p_i 
-\sum_{j=1}^h 2 q_j\right) \vert = \\
= \prod_{i=1}^n {{r_i + d_i} \choose r_i}(d_{n+1}+1) 
- (r_1 + \ldots + r_n + 2)(l + h),
\end{eqnarray*}
as expected.
\end{Proposition}
  
\proof Let $X := \bP^{r_1} \times \ldots \times \bP^{r_n} \times \bP^1$, 
$L := \mathcal{O}_X(d_1, \ldots, d_n, d_{n+1})$ and $Z := \cup_{i=1}^l 2 p_i 
\cup \cup_{j=1}^h 2 q_j$; we have to show that $H^1(L \otimes 
\mathcal{I}_{Z,X})=0$. Recall the natural exact sequence
$$
0 \to \mathcal{I}_{\mathrm{Res}_D(Z),X} \otimes L(-D) \to 
\mathcal{I}_{Z,Y} \otimes L \to \mathcal{I}_{Z \cap D, D} \otimes 
(L_{\vert D}) \to 0, 
$$
where $\mathrm{Res}_D(Z) = \cup_{i=1}^l 2 p_i \cup \cup_{j=1}^h  q_j$
denotes the residual scheme of $Z$ with respect to $D$; we are going 
to apply the usual Horace method by checking that 
\begin{eqnarray}
\label{first}
H^1 \left(\mathcal{I}_{\mathrm{Res}_D(Z),X} \otimes L(-D)\right) &=& 0 \\
\label{second}
H^1 \left(\mathcal{I}_{Z \cap D, D} \otimes (L_{\vert D})\right) &=& 0.
\end{eqnarray}
Since 
$$
\mathcal{I}_{\mathrm{Res}_D(Z),X}\otimes L(-D) = 
\mathcal{O}_{\bP^{r_1} \times \ldots \times 
\bP^{r_n} \times \bP^1}(d_1, \ldots, d_n, d_{n+1}-1) \left(-\sum_{i=1}^l 2 p_i 
-\sum_{j=1}^h  q_j\right),
$$
in order to get (\ref{first}) it is enough to see that the $h$ points $q_j$ 
on $D$ impose independent conditions on the non-special linear system 
$\vert \mathcal{O}_{\bP^{r_1} \times \ldots \times \bP^{r_n}\times \bP^1}
(d_1, \ldots,$ $d_n, d_{n+1}-1) \left(-\sum_{i=1}^l 2 p_i\right) \vert$. 
Indeed, if this were not the case, 
every divisor in such a linear system passing through 
$q_1, \ldots, q_{h-1}$ should contain $D$, hence 
\begin{eqnarray*}
\dim \vert \mathcal{O}_{\bP^{r_1} \times \ldots \times \bP^{r_n} \times \bP^1}
(d_1, \ldots, d_n, d_{n+1}-2)\left(-\sum_{i=1}^l 2 p_i\right) \vert \ge \\
\ge \dim \vert \mathcal{O}_{\bP^{r_1} \times \ldots \times 
\bP^{r_n}}(d_1, \ldots, d_n-1) \left(-\sum_{i=1}^l 2 p_i 
-\sum_{j=1}^{h-1} q_j\right) \vert \ge \\
\ge \prod_{i=1}^n {{r_i + d_i} \choose r_i}(d_{n+1}) 
- (r_1 + \ldots + r_n + 2)l - (h-1),
\end{eqnarray*}
in contradiction with our assumptions. Finally, notice that 
$$
\mathcal{I}_{Z \cap D, D} \otimes (L_{\vert D}) = 
\mathcal{O}_{\bP^{r_1} \times \ldots \times \bP^{r_n}}
(d_1, \ldots, d_n) \left(-\sum_{j=1}^h 2 q_j\right),
$$
so (\ref{second}) follows from our assumptions too. 

\qed

\begin{Proposition}\label{weakly}
Fix integers $n \ge 1$, $r_1, \ldots, r_n \ge 1$, $d_1, \ldots, d_n \ge 2$,
$1 \le s \le \left[ \frac{ \prod_{i=1}^{n} {{r_i + d_i} \choose r_i} 
(d_{n+1}+1)}
{(r_1 + \ldots + r_n + 2)} \right]$, $h_0 := \left[ \frac{ \prod_{i=1}^n 
{{r_i + d_i} \choose r_i}}{(r_1 + \ldots + r_n+1)} \right]$, 
$t_0$ such that $1 \le s - t_0 h_0 \le h_0$, and $d_{n+1} \ge t_0 + 3$.
Assume that 
\begin{eqnarray*}
\dim \vert \mathcal{O}_{\bP^{r_1} \times \ldots \times \bP^{r_n} \times \bP^1}
(d_1, \ldots, d_n, d_{n+1}-t+1) \left(-\sum_{i=1}^{s- t h_0} 2 p_i
-\sum_{j=1}^{h_0} 2 q_j\right) \vert = \\
= \prod_{i=1}^n {{r_i + d_i} \choose r_i}(d_{n+1}-t+2) 
- (r_1 + \ldots + r_n + 2)(s-(t-1)h_0) 
\end{eqnarray*}
for every $1 \le t \le t_0$, and  
\begin{eqnarray*}
\dim \vert \mathcal{O}_{\bP^{r_1} \times \ldots \times 
\bP^{r_n} \times \bP^1}(d_1, \ldots, d_n, d_{n+1}-t_0) \left(- 2 p_1
-\sum_{j=1}^{s - t_0 h_0 -1} 2 q_j\right) \vert = \\
= \prod_{i=1}^n {{r_i + d_i} \choose r_i}(d_{n+1}-t_0+1) 
- (r_1 + \ldots + r_n + 2)(s-t_0h_0),
\end{eqnarray*}
where the $p_i$'s are general points in 
$\bP^{r_1} \times \ldots \times \bP^{r_n} \times \bP^1$ 
and the $q_i$'s are general points on a divisor $D$ of type 
$(0, \ldots, 0, 1)$.
Then a general divisor in the linear system 
\begin{equation}\label{system}
\vert \mathcal{O}_{\bP^{r_1} \times \ldots \times \bP^{r_n} \times \bP^1}
(d_1, \ldots, d_n, d_{n+1}) \left(-\sum_{i=1}^s 2 p_i\right) \vert 
\end{equation} 
has only ordinary double points at $p_1, \ldots, p_s$  
and is elsewhere smooth. In particular, the Segre embedding 
of $V_{r_1,d_1} \times \ldots \times V_{r_n,d_n} \times V_{1, d_{n+1}}$ 
is not $(s-1)$-weakly defective. 
\end{Proposition}

\proof Since the linear system 
$$
\vert \mathcal{O}_{\bP^{r_1} \times \ldots \times \bP^{r_n} \times \bP^1}
(d_1, \ldots, d_n, d_{n+1}) \left(-\sum_{i=1}^{s- h_0} 2 p_i
-\sum_{j=1}^{h_0} 2 q_j\right) \vert
$$
has the expected dimension, a general divisor in (\ref{system}) 
specializes to $E + D$, where $E$ is a general divisor in 
$$
\vert \mathcal{O}_{\bP^{r_1} \times \ldots \times \bP^{r_n} \times \bP^1}
(d_1, \ldots, d_n, d_{n+1}-1) \left(-\sum_{i=1}^{s- h_0} 2 p_i \right)
\vert.
$$
Hence by \cite{ChiCil:02}, Theorem~1.4, we are reduced to prove that 
$E$ has only ordinary double points. Moreover, by repeating $t_0$ 
times exactly the same argument, we are reduced to prove that the 
general divisor $F$ in  
$$
\vert \mathcal{O}_{\bP^{r_1} \times \ldots \times \bP^{r_n} \times \bP^1}
(d_1, \ldots, d_n, d_{n+1}-t_0) \left(-\sum_{i=1}^{s- t_0 h_0} 2 p_i 
\right) \vert
$$
has only ordinary double points. Once again, let $F$ degenerate to 
$G+D$, where $G$ is a general divisor in 
$$
\vert \mathcal{O}_{\bP^{r_1} \times \ldots \times \bP^{r_n} \times \bP^1}
(d_1, \ldots, d_n, d_{n+1}-t_0-1) \left(- 2 p_1 \right) \vert.
$$
We claim that $G$ has an isolated ordinary double point in $p_1$.
Indeed, by \cite{ChiCil:02}, Theorem~1.4, it is enough to check that 
the Segre embedding of $V_{r_1,d_1} \times \ldots \times V_{r_n,d_n} 
\times V_{1, d_{n+1}-t_0-1}$ is not $0$-weakly defective, but this is clear, 
since by \cite{ChiCil:02}, Remark~3.1~(ii), $0$-weakly defective 
varieties contain many lines, here instead we have $d_i \ge 2$ for 
$i=1, \ldots, n-1$, and $d_{n+1}-t_0-1 \ge 2$. Hence the claim is established 
and we conclude by \cite{ChiCil:02}, Theorem~1.4. 

\qed  

\emph{Proof of Corollary~\ref{three}.} By \cite{CGG:03}, Theorem~2.1 and 
Theorem~2.5, if $d_i \ge 3$ for $i=1,2$, then both the linear systems
$\mathcal{O}_{\bP^1 \times \bP^1}(d_1,d_2)(-\sum 2 p_i)$ and 
$\mathcal{O}_{\bP^1 \times \bP^1 \times \bP^1}
(d_1,d_2,d_3)(-\sum 2 p_i)$ have always the expected dimension. 
In particular, $\Sigma_{1, \underline{d}}$ is never $k$-defective 
and for every integer $t$ such that $t \le d_3 - 3$ we have 
\begin{eqnarray*}
\dim \vert \mathcal{O}_{\bP^1 \times \bP^1 \times \bP^1} (d_1, d_2, d_3 -t -2) 
\left( - \sum_{i=1}^l 2 p_i \right) \vert \\
= (d_1+1)(d_2+1)(d_3-t+1)-4l \\
< (d_1+1)(d_2+1)(d_3-t+1)-4l - h
\end{eqnarray*}
for $h < (d_1+1)(d_2+1)$, so all assumptions of Proposition~\ref{horace} are 
satisfied. As a consequence, we can apply Proposition~\ref{weakly} 
(just notice that $t_0 \le d_3 - 3$ by (\ref{assumption})), and deduce that 
$\Sigma_{1, \underline{d}}$ is not $(k-1)$-defective. 
Now the claim follows from Theorem~\ref{main}. 

\qed

\vspace{0.5cm}

\noindent
Claudio Fontanari \newline
Universit\`a degli Studi di Trento \newline
Dipartimento di Matematica \newline
Via Sommarive 14 \newline
38050 Povo (Trento) \newline
Italy \newline
e-mail: fontanar@science.unitn.it

\end{document}